\mag=\magstep1
\documentstyle{amsppt}
\input amsppt1
\pageheight{23true cm}
\pagewidth{17true cm}

\parindent=4mm
\parskip=3pt plus1pt minus.5pt
\nologo\NoRunningHeads\NoBlackBoxes

\def\si{\roman{sin}}
\def\co{\roman{cos}}
\def\i{\looparrowright}
\def\e{\hookrightarrow}

\def\pa{ pair of $n$-knots}
\def\f{\flushpar }
\def\nl{\newline }
\def\np{\newpage }
\def\x{\times }
\def\p{\bf Proof of  }

Intersectional pairs of $n$-knots, \newline
local moves of $n$-knots, and \newline
their associated invariants of $n$-knots  

\vskip1cm

Eiji Ogasa

\vskip1cm
pqr100pqr100\@yahoo.co.jp

\vskip1cm
Abstract.
Let $n$ be an integer$>0$. 
Let $S^{n+2}_1$ (respectively, $S^{n+2}_2$) be the $(n+2)$-sphere 
embedded in the$(n+4)$-sphere  $S^{n+4}$. 
Let $S^{n+2}_1$ and $S^{n+2}_2$ intersect transversely. 
Suppose that the smooth submanifold $S^{n+2}_1\cap S^{n+2}_2$ in $S^{n+2}_i$ 
is PL homeomophic to the $n$-sphere. 
Then 
$S^{n+2}_1\cap S^{n+2}_2$ in $S^{n+2}_i$ 
is an $n$-knot $K_i$. 
We say that the pair $(K_1,K_2)$ of $n$-knots is realizable. 

We consider the following problem in this paper. 
Let $A_1$ and $A_2$ be $n$-knots. 
Is the pair $(A_1,A_2)$ of $n$-knots realizable?

We give a complete characterization.

\vskip1cm
Chapter I is accepted in \newline
Mathematical Research Letters, 1998, 5, 577-582.\newline
This manuscript is not the published version.

Chapter I is a summary of Chapter II. 
Chapter II is beased on the author's PhD thesis (University of Tokyo, 1996)
University of Tokyo preprint series UTMS 95-50
 is a preprint of Chapter II.

\newpage
\vskip5cm

Chapter I

\vskip5cm

\topmatter
\title
Intersectional pairs of $n$-knots, \\
local moves of $n$-knots, and \\
their associated invariants of $n$-knots  
\endtitle
\affil
Eiji Ogasa\\
ogasa\@ms.u-tokyo.ac.jp, 
ogasa\@max.math.brandeis.edu
\endaffil
\thanks{This research was partially supported by Research Fellowships 
of the Promotion of Science for Young Scientists.\newline}
\endtopmatter

\document

\baselineskip11pt

\head 1. Introduction \endhead     

Our first purpose is to discuss the following problem.

Let $S^{n+2}_1$ and $S^{n+2}_2$ be $(n+2)$-spheres 
embedded in the  $(n+4)$-sphere $S^{n+4}$ 
($n\geqq1$) 
which  intersect transversely. 
If we assume $M=S^{n+2}_1\cap S^{n+2}_2$ is PL homeomorphic to the single standard $n$-sphere, 
we obtain a pair of  $n$-knots, 
$M$ in $S^{n+2}_1$ and $M$ in $S^{n+2}_2$.   
We consider which pairs of $n$-knots we obtain as above. 
That is, let $(K_1,K_2)$ be a pair of $n$-knots. 
Then we consider whether the pair of $n$-knots ($K_1$, $K_2$) 
is obtained as above.  

We give a complete answer to this problem (Theorem 3.1).

In order to get the complete answer, we introduce a local move of $n$-knots ($n\geqq1$). 
Furthermore, we show a relation between the local move and 
some invariants of $n$-knots (Theorem 4.1 and Corollary 4.2). 

Our second purpose is to discuss the relation between 
the local move and the invariants of $n$-knots.  
In the case of 1-links, 
there is a great deal known about 
relations between local moves and knot invariants. 
(See e.g. [V][Wi][Ka2].)
Our discussion is a high dimensional version of this theory.

\head 2. Definitions \endhead

An {\it (oriented) (ordered) $m$-component n-(dimensional) link}
 is a smooth, oriented submanifold $L=\{K_1,...,K_m\}$ of $S^{n+2}$, 
  which is the ordered disjoint union of 
$m$ connected oriented submanifolds, each PL homeomorphic 
to the standard $n$-sphere. If $m=1$, then $L$ is called a {\it knot}.
(This definition is used often. See e.g. [CO],[L1],[L3].)

Let $L_1$ and $L_2$ be $n$-links.
$L_1$ is said to be equivalent to $L_2$ if there exists an orientation 
preserving diffeomorphism $h$ of $S^{n+2}$ such that
$h\vert L_1$ is an orientation preserving diffeomorphism from $L_1$ to $L_2$.

We work in the smooth category. 

\f{\bf Definition}
$(K_1, K_2)$ is called a {\it pair of $n$-knots }
if $K_1$ and $K_2$ are $n$-knots.  
$(K_1, K_2, X_1, X_2)$ is called a {\it 4-tuple of $n$-knots 
and $(n+2)$-knots} or a {\it 4-tuple of $(n,n+2)$-knots} 
if 
($K_1, K_2$) is a pair of $n$-knots 
and 
$ X_1$ and $X_2$ are $(n+2)$-knots 
diffeomorphic to the standard $(n+2)$-sphere. 
($n\geqq1$).

\definition  {Definition}
A 4-tuple of $(n,n+2)$-knots $(K_1, K_2, X_1, X_2)$   
is said to be  {\it realizable } if there exists a smooth transverse 
immersion  $f:S^{n+2}_1\coprod S^{n+2}_2 \looparrowright S^{n+4}$  
satisfying the  following conditions. 
($n\geqq1$). 
\roster
\item
The intersection  $\Sigma=f(S^{n+2}_1)\cap f(S^{n+2}_2)$ is 
PL homeomorphic to the standard $n$-sphere. 
\item  
$f^{-1}(\Sigma)$ in $S^{n+2}_i$ defines an $n$-knot $K_i$ ($i=1,2$).
\item
 $f\vert S^{n+2}_i$ is an embedding. 
 $f(S^{n+2}_i)$ in $S^{n+4}$ is equivalent to $X_i$ (i=1,2).
\endroster   

A pair of $n$-knots $(K_1, K_2)$ is said to be {\it realizable } 
or is called an {\it intersectional pair of $n$-knots }  
if there is a realizable 4-tuple of $(n,n+2)$-knots $(K_1, K_2, X_1, X_2)$.   
\enddefinition

\head 3. Intersectional pair of $n$-knots \endhead

Our main theorem is:

\proclaim{Theorem 3.1}
A \pa       \hskip1mm $(K_1, K_2)$ 
($n\geqq1$)  
is realizable  if and only if 
\nl $(K_1, K_2)$  satisfies the condition that 

 \hskip15mm
$\cases
\text{$(K_1, K_2)$ is arbitrary} & \text{if $n$ is even,} \\  
\text{Arf($K_1$)=Arf($K_2$)} & \text{if $n=4m+1$,}\\
\text{$\sigma(K_1$)=$\sigma(K_2$)} & \text{if $n=4m+3$,}\\
\endcases$
($m\geqq 0$). 
\endproclaim

There is a mod 4 periodicity in dimension. 
It is similar to the periodicity 
in knot cobordism theory ([L1]) and 
surgery theory (See e.g. [Br][Wa][CS][We]).  

We have the following result on the realization of 4-tuples of ($n,n+2$)-knots. 

\proclaim{Theorem 3.2}
A 4-tuple of ($n,n+2$)-knots  $(K_1, K_2, X_1,X_2)$ is realizable  
if $K_1$  and  $K_2$ are slice.  
($n\geqq1$). 
In particular, 
if $n$ is even, an arbitrary 4-tuple of  ($n,n+2$)-knots 
 $(K_1, K_2, X_1, X_2)$ is realizable.
\endproclaim

{\bf Note.} 
(1) Kervaire proved  that all even dimensional knots are  slice ([Ke]).

(2) In [O1] the author discussed the case of two 3-spheres in a 5-sphere.
In [O2] the author discussed the case of the intersection of three 4-spheres in a 6-sphere.

\definition{Problem }
Which 4-tuples of $(2n+1,2n+3)$-knots are realizable ($n\geqq1$)?
\enddefinition

\head 4. High-dimensional pass-moves \endhead

In order to prove Theorem 3.1, we introduce a new local move for  high dimensional knots, 
the {\it high dimensional pass-move}. 
Pass-moves for 1-knots are discussed in p.146 of [Ka].

We define high dimensional pass-moves for $(2k+1)$-knots $\subset S^{2k+3}$.

\f{\bf Definition}
Take a trivially embedded 
$(2k+3)$-ball 
$B=B^{2k+2}\x[-1,1]$ 
in $S^{2k+3}$.
We define $J_+, J_-\subset B$ as follows. 
Refer to Figure 4.1. 

In  $\partial B^{2k+2}\x\{0\}$, 
take trivially embedded $S_1^{k}$, $S_2^{k}$  
such that lk$(S_1^{k}, S_2^{k})=1$.  
Let $N(S_*^{k})$ be a tubular neighborhood of $S_*^{k}$ 
in  $\partial B^{2k+2}\x\{0\}$.

Let $h^{k+1}$ be an $(2k+2)$-dimensional $(k+1)$-handle 
which is attached to $\partial B^{2k+2}\x\{0\}$ along $N(S_1^{k})$   
with the trivial framing  
and 
which is embedded trivially in $B^{2k+2}\x\{0\}$.     

Let $h^{k+1}_+$ (resp. $h^{k+1}_-$)
 be an $(2k+2)$-dimensional $(k+1)$-handle 
which is embedded 
in  $B=B^{2k+2}\x[0,1]$ (resp. $B=B^{2k+2}\x[-1,0]$)
and 
which is attached to $\partial B^{2k+2}\x\{0\}$ along $N(S_2^{k})$   
with the trivial framing.

Let $h^{k+1}_+\cap h^{k+1}_-=N(S_2^{k})$.  
Let $h^{k+1}_+\cap h^{k+1}$=$h^{k+1}_-\cap h^{k+1}$=$\phi$.

Let $J_+$ be a submanifold 
$\overline{ (\partial h^{k+1})-N(S_1^{k})} \amalg 
\overline{ (\partial h^{k+1}_+)-N(S_2^{k})}$
in  $B$. 

Let $J_-$ be a submanifold 
$\overline {(\partial h^{k+1})-N(S_1^{k})}\amalg 
\overline  {(\partial h^{k+1}_-)-N(S_2^{k})}$
in  $B$.

In Figure 4.1, 
we draw  $B=B^{2k+2}\x[-1,1]$  by using the projection to $B^{2k+2}\x\{0\}$.  


\vskip1cm
Figure 4.1.

You can obtain this figure 
by clicking `PostScript' in the right side of 
the cite of the abstract of this paper in arXiv 
(https://arxiv.org/abs/the number of this paper). 

You can also obtain it from the author's website,  
which can be found by typing his name in search engine. 

\vskip1cm

Let $K_+$, $K_-$  be $(2k+1)$-knots  $\subset S^{2k+3}$. 
We say that 
$K_+$ is obtained from $K_-$ by one {\it high dimensional pass-move} 
if there is a trivially embedded $(2k+2)$-ball $B\subset S^{2k+3}$ 
such that
$K_+\cap B$ is $J_+$ and $K_-\cap B$ is $J_-$.  

Let $K$, $K'$  be $(2k+1)$-knots  $\subset S^{2k+3}$. 
We say that 
$K$ is {\it pass-move equivalent} to $K'$ 
if there are $(2k+1)$-knots $K_1$,...,$K_\mu$ ($\mu\in\Bbb N$)
such that $K_i$ is pass-move equivalent to $K_{i+1}$.

We prove:
\proclaim{Theorem 4.1} 
For $(2k+1)$-knots $K_1$ and $K_2$, 
the following two conditions are equivalent. ($k\geqq1$.)
\roster
\item
There exists a $(2k+1)$-knot $K_3$ 
which is pass-move equivalent to $K_1$ and 
cobordant to $K_2$. 
\item
$K_1$ and $K_2$ satisfy the condition that 
$\cases
\roman{Arf}(K_1)=\roman{Arf}(K_2)
 & \text{when $k$ is even}\\ 
 \text{ 
$\sigma (K_1)\,=\,\sigma (K_2)$ 
 }
 & \text{when $k$ is odd.} 
\endcases$
\endroster
\endproclaim

The $k=0$ case of Theorem 4.1 follows from [Ka].

\proclaim{Corollary 4.2} 
Let $K_1$ and $K_2$ be $(2k+1)$-knots ($k\geqq1$.). 
Suppose that 
$K_1$ is pass-move equivalent to $K_2$. 
Then 
$K_1$ and $K_2$ satisfy the condition that 

$\cases
\roman{Arf}(K_1)=\roman{Arf}(K_2)
 & \text{when $k$ is even}\\ 
 \text{ 
$\sigma (K_1)=\sigma (K_2)$ 
 }
 & \text{when $k$ is odd.} 
\endcases$
\endproclaim 

{\bf Note.} In [O3] the author proved a relation between 
another local move of 2-knots and other invariants of 2-knots.

\head 5. Proof of Theorem 3.1 \endhead

We prove the following lemmas by explicit construction. 

\proclaim{Lemma 5.1}
Let $K$ be an $n$-knot. 
Then the pair of $n$-knots $(K,K)$ is realizable ($n\geqq1$).  
\endproclaim

\proclaim{Lemma 5.2}
Let $K_1$ and $K_2$ be $(2k+1)$-knots. 
Suppose that $K_1$ is pass-move equivalent to $K_2$. 
Then the pair of $(2k+1)$-knots $(K_1, K_2)$ is realizable ($k\geqq0$).  
\endproclaim

\proclaim{Lemma 5.3}
Let $K_1$, $K_2$ and $K_3$ be $n$-knots ($n\geqq1$). 
Suppose that the pair of $n$-knots $(K_1,K_2)$ is realizable 
and that $K_2$ is cobordant to $K_3$.  
Then the pair of $n$-knots $(K_1, K_3)$ is realizable.  
\endproclaim

Theorem 3.1 is deduced from Theorem 4.1 and Lemmas 5.1, 5.2, 5.3.

\head 6. Proof of Theorem 3.2 \endhead

It suffices to prove that 
a 4-tuple of ($n,n+2$)-knots $(K_1,K_2,T,T)$ is realizable, 
where $K_1$ is a slice $n$-knot, $K_2$ is the trivial $n$-knot, 
$T$ is the trivial $(n+2)$-knot. 

Any 1-twist spun knot is unknotted ([Z]).
Theorem 3.2 follows from this fact.

\head 7. The proof of Theorem 4.1\endhead

Every $p$-knot $(p>1)$ is cobordant to a simple knot. 
(See [L1] for a proof and the definition of simple knots. )
By using this fact, we prove that 
the $k\geqq1$ case of Theorem 4.1 can be deduced from Theorem 7.1.  

\proclaim{Proposition 7.1} 
For simple $(2k+1)$-knots $K_1$ and $K_2$, 
the following two conditions are equivalent. ($k\geqq1$.)
\roster
\item
$K_1$ is pass-move equivalent to $K_2$. 
\item
$K_1$ and $K_2$ satisfy the condition that 
$\cases
\roman{Arf}(K_1)=\roman{Arf}(K_2)
 & \text{when $k$ is even}\\ 
 \text{ 
$\sigma (K_1)=\sigma (K_2)$ 
 }
 & \text{when $k$ is odd.} 
\endcases$
\endroster
\endproclaim 

{\bf Proof of Proposition 7.1. }
(2)$\Rightarrow$(1).
$K_1$ bounds a Seifert hypersurface $V_1$ with a handle decomposition 
(one 0-handle)$\cup$((k+1)-handles). 
Take a Seifert matrix associated with $V_1$. 
By using high dimensional pass moves, we can change the Seifert matrix 
without changing the diffeomorphism type of $V_1$. 
Thus we obtain a $(2k+1)$-knot $K'_2$ whose Seifert matrix is same as 
the Seifert matrix of $K_2$ 
if (2) holds.  
By the classification theorem of simple knots by [L2], 
$K'_2$ is equivalent to $K_2$. 

(1)$\Rightarrow$(2). 
Suppose that ($2k+1$)-knots $K_*\subset S_*^{2k+3}$ bounds 
a Seifert hyper surface $V_*$. Note $V_*$ are ($2k+2$)-manifolds.
There is a compact oriented parallelizable $(2k+4)$-manifold $P$ 
whose boundary is $S_1^{4k+3}\amalg S_2^{4k+3}$ 
containning compact oriented  $(2k+3)$-manifold $Q$ 
whose boundary is $V_1\cup (S^{2k+1}\times [1,2]) \cup V_2$. 
( Here, $\partial V_*$ is $K_*$ and $S^{2k+1}\times\{*\}$ is $K_*$. ) 
We use characteristic classes and intersection products to prove (1)$\Rightarrow$(2).

\head 8. Intersectional pair of submanifolds \endhead

In \S1 suppose $M$ is not PL homeomorphic to the standard sphere.
Then we obtain a pair of submanifolds, $M$ in $S^{n+2}_i$ ($i=1,2$).   

Let $N$ be a closed oriented manifold.
$(K_1, K_2)$ is called {\it a pair of submanifolds (diffeomorphic to $N$)} 
if $K_i$ is a submanifold of $S^{n+2}$ 
diffeomorphic to $N$.

Let $(K_1, K_2)$ be 
a pair of submanifolds diffeomorphic to $M$.  
We say $(K_1, K_2)$ is an {\it intersectional pair } 
if the submanifold $K_i$ is equivalent to the submanifold 
$M=S^{n+2}_1\cap S^{n+2}_2$ in $S^{n+2}_i$ 
as in \S1 ($i=1,2$).

It is natural to ask the following problem.  

\definition{Problem 8.1}   
Which pairs of submanifolds are intersectional pairs?
\enddefinition

The author can prove the following results. 

When $n$ is even, not all pair of submanifolds as above are realizable. 

When $n=4m+3$, we can define the signature as in the knot case and the signature is 
an obstruction. Therefore not all pairs are realizable. 
When $n=3$, 
$(K_1,K_2)$ is realizable if and only if $\sigma(K_1)$$=\sigma(K_2)$. 
When $n\neq3$, 
$\sigma(K_1)$$=\sigma(K_2)$ does not imply $(K_1,K_2)$ is realizable in general. 

When $n=4m+1$, there is a closed oriented manifold $M$  such that 
if $K_1$ and $K_2$ are PL homeomorphic to $M$, then ($K_1,K_2$) is realizable. 
In other words, 
there is no invariant corresponding to the Arf invariant as in the knot case. 
Of course, not all pairs are realizable.



\newpage

\vskip5cm

 Chapter II

\vskip5cm

\topmatter
\title
On the intersection of spheres in a sphere II:\\High dimensional case
\endtitle
\author
Eiji Ogasa
\endauthor
\thanks{This research was partially suppported by Reseach Fellowships 
of thePromotion of Science for Young Scientists.\newline}

\endtopmatter

\document
\f{\bf Abstract.}
Consider transverse immersions $f:S^{n+2}_1\amalg S^{n+2}_2\i S^{n+4}$ 
such that  $f\vert S^{n+2}_i$ is an embedding and 
the intersection $f(S^{n+2}_1)$$\cap f(S^{n+2}_2)$ 
is PL homeomorphic to the standard $n$-sphere. 
($n\geqq1$). 
Then we obtain a pair of $n$-knots, 
$f^{-1}$($f(S^{n+2}_1)$$\cap f(S^{n+2}_2)$) in  $S^{n+2}_i$ ($i=1,2$).
We determine which pair of $n$-knots are obtained as above.  
Roughly speaking, our result is characterized 
by the Arf invariant and the signature. 
We find  a mod 4 periodicity in the dimension $n$.

\baselineskip11.7pt
\head 1.Introduction and Main results \endhead

Let $S^{n+2}_1$ and $S^{n+2}_2$ be the $(n+2)$-spheres 
embedded in the  $(n+4)$-sphere $S^{n+4}$ 
($n\geqq1$) 
and  intersect transversely. 
Here, the orientation of  the intersection $M$ is induced naturally.
If we assume $M$ is PL homeomorphic to the single standard $n$-sphere, 
we obtain a pair of  $n$-knots, $M$ in $S^{n+2}_i$ ($i=1,2$),  
and a pair of  $(n+2)$-knots, $S^{n+2}_i$ in  $S^{n+4}$  ($i=1,2$).

Conversely, let ($K_1$, $K_2$) be a pair of $n$-knots.  
It is natural to ask whether ($K_1$, $K_2$) is obtained as above.  
In this paper we give a complete answer to the above question.
Furthermore we discuss somewhat 
which 4-tuple($K_1$, $K_2$, $S^{n+2}_1$, $S^{n+2}_2$) are realizable.

To state our results we need some definitions.

An {\it (oriented) (ordered) $m$-component n-(dimensional) link}
 is a smooth, oriented submanifold $L=\{K_1,...,K_m\}$ of $S^{n+2}$, 
  which is the ordered disjoint union of $m$ manifolds, each PL homeomorphic 
to the standard $n$-sphere. (If $m=1$, then $L$ is called a {\it knot}.)  
We say that $m$-component $n$-dimensional links, 
$L_0$ and $L_1$,  are said to be 
{\it (link-)concordant} or {\it  (link-)cobordant} 
if there is a smooth oriented submanifold  
$\widetilde{C}$=\{$C_1$ ,...,$C_m$\} of $S^{n+2}\times [0,1]$,  
which meets the boundary transversely in $\partial \widetilde{C}$, 
      is PL homeomorphic to  $L_0 \times[0,1]$ 
  and meets $S^{n+2}\times \{l\}$ in $L_l$ ($l=0,1$).
 (See [CO]).

\f{\bf Definition 1.1.}
$(K_1, K_2)$ is called a {\it pair of $n$-knots }
if $K_1$ and $K_2$ are $n$-knots.  
$(K_1, K_2, X_1, X_2)$ is called a {\it 4-tuple of $n$-knots 
and $(n+2)$-knots} or a {\it 4-tuple of $(n,n+2)$-knots} 
if $K_1$ and $K_2$ compose a pair of $n$-knots $(K_1, K_2)$ 
and 
$ X_1$ and $X_2$ are $(n+2)$-knots 
diffeomorphic to the standard $(n+2)$-sphere. 
($n\geqq1$).

\definition  {Definition 1.2}
A 4-tuple of $(n,n+2)$-knots $(K_1, K_2, X_1, X_2)$   
is said to be  {\it realizable } if there exists a smooth transverse 
immersion  $f:S^{n+2}_1\coprod S^{n+2}_2 \looparrowright S^{n+4}$  
satisfying the  following conditions. 
($n\geqq1$). 
\roster
\item
 $f\vert S^{n+2}_i$  defines $X_i$ (i=1,2).
\item
The intersection  $\Sigma=f(S^{n+2}_1)\cap f(S^{n+2}_2)$ is 
PL homeomorphic to the standard $n$-sphere. 
\item 
$f^{-1}(\Sigma)$ in $S^{n+2}_i$ defines an $n$-knot $K_i$ ($i=1,2$).
\endroster
A pair of $n$-knots $(K_1, K_2)$ is said to be  {\it realizable } if 
there is a realizable 4-tuple of $(n,n+2)$-knots $(K_1, K_2, X_1, X_2)$.   
Then  $f$ is called an {\it immersion to realize} 
$(K_1, K_2, X_1, X_2)$ or $(K_1, K_2)$. 
\enddefinition
The following theorem characterizes the realizable pair of $n$-knots.

\proclaim{Theorem 1.3}
A \pa       \hskip1mm $(K_1, K_2)$ 
($n\geqq1$)  
is realizable  if and only if 
\nl $(K_1, K_2)$  satisfies the condition that 

 \hskip15mm
$\cases
\text{$(K_1, K_2)$ is arbitrary} & \text{if $n$ is even,} \\  
\text{Arf($K_1$)=Arf($K_2$)} & \text{if $n=4m+1$,}\\
\text{$\sigma(K_1$)=$\sigma(K_2$)} & \text{if $n=4m+3$,}\\
\endcases$
($m\geqq 0$,$m\in\Bbb Z$). 
\endproclaim

There exists a mod 4 periodicity in dimension. 
It is similar to the periodicity 
in the  knot cobordism theory and the surgery theory.

We have the following results on the realization of 4-tuple of ($n,n+2$)-knots. 

\proclaim{Theorem 1.4}
A 4-tuple of ($n,n+2$)-knots  $(K_1, K_2, X_1,X_2)$ is realizable  
if $K_1$  and  $K_2$ are slice.  
($n\geqq1$). 
\endproclaim
\f 
Kervaire proved  that all even dimensional knots are  slice ([K]).
Hence we have:
\proclaim{Corollary 1.5}
If $n$ is even, an arbitrary 4-tuple of  ($n,n+2$)-knots 
 $(K_1, K_2, X_1, X_2)$ is realizable.
\endproclaim

The author discussed related topics in [O1], [O2], and [O3].

This paper is organized as follows. 
In \S 2 
 we introduce a new knotting operation,  the high dimensional pass-move,  
 and  state  its relation to the Arf invariant and the signature.   
In \S 3 
we discuss a sufficient condition for the realization of 
pair of odd dimensional knots. 
In \S 4
 we discuss a necessary condition  
 for the realization of pair of ($4m+1$)-knots. 
In \S 5
 we discuss a necessary condition 
 for the realization of pair of ($4m+3$)-knots. 
In \S 6
we prove Theorem 1.4 which induces a necessary and sufficient condition  
 for the realization of pair of even dimensional knots. 
Theorem 1.3 follows from \S 2-6. 

The author would like to thank Professor Takashi Tsuboi for his advise.

\head 2.
High dimensional pass-moves 
\endhead

In this section we introduce a new knotting operation 
for high dimensional knots. 
The 1-dimensional case of Definition 2.1 
is disscussed in P. 146 of [Kf].

\definition{Definition 2.1}
Take a ($2k+1$)-knot $K$ $(k\geqq0)$. 
Let  $K$ be defined by a smooth embedding 
$g:\Sigma^{2k+1}\e S^{2k+3}$, 
where  $\Sigma^{2k+1}$ is PL homeomorphic to the standard $(2k+1)$-sphere. 
Let
$D^{k+1}_x$=$\{(x_1$$,...,x_{k+1})$$\vert$ $\Sigma x_i^2$$<1\}$ and 
$D^{k+1}_y$=$\{(y_1$$,...,y_{k+1})$$\vert$ $\Sigma y_i^2$$<1\}$.
Let
$D^{k+1}_x(r)$=$\{(x_1$ $,...,x_{k+1})$$\vert$ $\Sigma x_i^2$$\leqq r^2\}$ and 
$D^{k+1}_y(r)$=$\{(y_1$ $,...,y_{k+1})$$\vert$ $\Sigma y_i^2$$\leqq r^2\}$.
A local chart $(U,\phi)$ of $S^{2k+3}$ is called a 
{\it pass-move-chart of $K$} if it satisfies the following conditions. 
\roster
\item
 $ \phi(U)\cong$ $\Bbb R^{2k+3}$=
 $(0,1)\times D^{k+1}_x \times D^{k+1}_y$  
\item
 $ \phi( g(\Sigma^{2k+1}) \cap U)$  
  = $[\{\frac{1}{2}\}\times D^{k+1}_x \times \partial D^{k+1}_y(\frac{1}{3})]$  $\amalg$ 
  $[\{\frac{2}{3}\}\times\partial D^{k+1}_x(\frac{1}{3}) \times D^{k+1}_y]$  
\endroster 
Let $g_U:\Sigma^{2k+1}\e S^{2k+3}$ be an embedding such that: 
 \roster
\item
$g\vert\{\Sigma^{2k+1}-g^{-1}(U)\}=g_U\vert\{\Sigma^{2k+1}-g^{-1}(U)\}$, and 
\item
 $ \phi( g_U(\Sigma^{2k+1})\cap U)$
= $[\{\frac{1}{2}\}\times D^{k+1}_x \times\partial D^{k+1}_y(\frac{1}{3})]$  
$\amalg$ 

 [$\{\frac{2}{3}\}  \times 
 \partial D^{k+1}_x(\frac{1}{3})      \times $
 $(D^{k+1}_y- D^{k+1}_y(\frac{1}{2}))$]  

$\cup$ 
[$[\frac{1}{3},\frac{2}{3}]\times\partial D^{k+1}_x(\frac{1}{3}) 
\times \partial D^{k+1}_y(\frac{1}{2})$]

$\cup$ 
$[\{\frac{1}{3}\}\times\partial D^{k+1}_x(\frac{1}{3}) 
\times D^{k+1}_y(\frac{1}{2})\}]$  
\endroster
Let $K_U$ be the  $(2k+1)$-knot defined by   $g_U.$
Then we say that $K_U$ is obtained from $K$ by the 
{\it (high dimensional) pass-move in $U$.}
We say that 
$(2k+1)$-knot $K$ and $K'$ are 
 {\it (high dimensional) pass-move equivalent} 
  if there  exist $(2k+1)$-knots 
$K$=$K_1$,$K_2$...,$K_q$,$K_{q+1}$=$K'$ 
and 
$K_{i+1}$ is obtained from $K_i$ 
by the high dimensional pass-move 
in a pass-move-chart of $K_i$ $(i=1,...,q)$. 


\enddefinition

High dimensional pass-moves  have the follwing relation 
with the Arf invariant and the signature of knots. 
 (The case of $k=0$ follows from [Kf].)
We prove:

\proclaim{Theorem 2.2}
For $(2k+1)$-knots $K_1$ and $K_2$, 
the following two conditions are equivalent. ($k\geqq0$.)
\roster

\item
There exists a $(2k+1)$-knot $K_3$ 
which is pass-move equivalent to $K_1$ and 
cobordant to $K_2$. 

\item
$K_1$ and $K_2$ satisfy the condition 
$\cases
\text{ 
 Arf$(K_1)$=Arf$(K_2)$ 
 }
 & \text{when $k$ is even}\\ 
 \text{ 
$\sigma (K_1)=\sigma (K_2)$ 
 }
 & \text{when $k$ is odd.} 
\endcases$
\endroster
\endproclaim


Organization of the proof of Theorem 2.2 is as follows.
Obviously Theorem 2.2 is equivalent to the following Claim 2.2.1 and  2.2.2.

\f
{\bf Claim 2.2.1.}
 {\sl If (2) of Theorem 2.2 holds, then (1) of Theorem 2.2 holds.}

\f
{\bf Claim 2.2.2.}
 {\sl If (1) of Theorem 2.2 holds, then (2) of Theorem 2.2 holds. }

\f 
In this section we prove Claim 2.2.1. (We use Claim 2.2.1 in \S 3.)
We use the results of \S 3, 4 and 5 and  prove Claim 2.2.2. 
(Note. We don't use Claim 2.2.2 in the proof of \S 3, 4 and 5.)
The proof of Claim 2.2.2 is written in  
\S 5.A. after \S 5.


We begin the proof of  Claim 2.2.1.
We need the following Lemma 2.3. 
We prove:

\proclaim{Lemma 2.3}
If a ($2k+1$)-knot $K$ ($k\geqq0$) satisfy the condition 
\nl
$(*)$ 
$\cases
\text{ 
 Arf($K)$=0} 
 & \text{when $k$ is even} 
 \\ 
 \text{ 
$\sigma (K)=0$ 
 }
 & \text{when $k$ is odd,} 
\endcases$
then there exists a ($2k+1$)-knot $\widetilde{K}$ 
which is pass-move equivalent to the trivial knot 
and cobordant to the ($2k+1$)-knot $K$. 
\endproclaim

Before proving Lemma 2.3, 
We prove: 

\f{\bf Claim.}
 Lemma 2.3 induces  
 Claim 2.2.1.

 Proof. 
 $(-K^*_1)\sharp K_2$ satisfies the condition $(*)$.  
By Lemma 2.3, there exists a ($2k+1$) knot $\widetilde{K}$ 
which is pass-move equivalent to the trivial knot 
and cobordant to  $(-K^*_1)\sharp K_2$
Define $K_3$  to be $K_1\sharp$ $\widetilde{K}$. 
Then the following (1) and (2) hold.
(1)$K_3$ is pass-move equivalent to $K_1$. 
(2)$K_3$ = $K_1\sharp$$\widetilde{K}$ is cobordant to 
 $K_1\sharp (-K^*_1)\sharp K_2$ and to $K_2$. \qed

Before proving Lemma 2.3, we review some definitions.
(See  [L2] and [Kw] for detail.)
We first review on the definition of  the Seifert matrix. 
Let $K$ be a ($2k+1$)-knot and $F$ a Seifert hypersurface. 
Let $F \times[-1,1]$ be embedded in  $S^{2k+3}$ 
so that $F \times\{0\}$ coincides with $F$ and 
 the standard orientation of [-1,1] coincides with the orientation of 
the normal bundle induced from that of $F$ and that of $S^{2k+3}$. 
For ($k+1$)-cycles $u$ and $v$ in $F$, define $\theta(u,v)$  
to be 
$lk(u,v\times\{1\})$ 
in $S^{2k+3}$. 
Let $z_1,...,z_{p}$ be $(k+1)$-cycles in $F$ which represent  basis of 
$H_{k+1}(F:\Bbb Z)$/(Torsion part). 
Define the {\it Seifert matrix $A=\{a_{ij}\}$ of $K$ associated 
with $F$ and $z_i$} to be $a_{ij}$ =$\theta (z_i,z_j)$. 
Here, recall the following lemma 2.4.

\proclaim{Lemma 2.4}(Well-known.) 
\roster
\item 
Let $K$ be a $(2k+1)$-knot ($k\geqq 0$)
with a Seifert hypersurface $F$. 
Let  $u$ and $v$ be ($k+1$)-cycles in $F$ $\subset S^{2k+3}$.
We have: 

\hskip15mm$ \theta(u,v)+(-1)^{k+1}\theta(v,u)=u\cdot v,$

\f where $u\cdot v$ is the intersection number in $F$.

\item
For vanishing ($k+1$)-cycles $\mu$ and $\nu$ 
 in $S^{2k+3}$ ($k\geqq0$), 
 where  $\mu\cap \nu=\phi$,   

\hskip15mm$lk(\mu,\nu)=(-1)^{k}lk(\nu,\mu). $

\item
Let $K$ be a $(2k+1)$-knot with a Seifert matrix $A$  ($k \geqq0$). 
For appropriate Seifert surfaces, 
the Seifert matrixes of $-K,$ that of $K^*$ and that of  $-K^*$ are   
$(-1)^{k}$$\{^t\hskip-1mm A\}$, $(-1)^{(k+1)}$$\{^t\hskip-1mm A\}$ and $-A$,
 respectively.
\endroster
\endproclaim

Recall the following theorem. 
\proclaim{Theorem 2.5}([L1])
(2k+1)-knots ($k\geqq1$, $k\neq0$) 
$K_1$ and $K_2$ are cobordant if and only if 
for a Seifert matrix $A_i$ of $K_i (i=1,2)$, 
$\pmatrix
A_1&O\\
O&{-A_2}\\
\endpmatrix$
is congruent to 
$\pmatrix
O&N_1\\
N_2&N_3\\
\endpmatrix,$  
where $N_i$ are same size. 
\endproclaim

We next review on the definition of the Arf invariant.
Let $K$ be a ($4m+1$)-knot. 
For $(2m+1)$-cycles $x$ in $F$,
 define $q(x)\in\Bbb Z_2 $ to be  $\theta(x,x)$ mod 2. 
Let $x_1,...,x_p$, $y_1,...,y_p$ be  
$(2m+1)$-cycles in $F$ which represent sympletic basis of 
$H_{2m+1}(F:\Bbb Z)$/(Torsion part), i.e., 
basis such that 
(1) $x_i\cdot y_i=1$ for all $i$,  
(2)$x_i\cdot x_j=0$, $y_i\cdot y_j=0$ for all $(i,j)$, and 
(3)$x_i\cdot y_j=0$ for $i\neq j$. 
Note that symplectic basis always exist.  
Then define $\roman{Arf}(K)=\sum_{i=1}^p q(x_i)q(y_i) \in \Bbb Z_2.$

At the end we review on the definition of the signature of knots.
Define {\it the signature} 
$\sigma(K)$ of $K$ to be the signature of $A+{^t\hskip-1mm}A$.
Recall that,  when $k$ is odd, $\sigma(K)$= $\sigma(F)$.

We now begin the proof of Lemma 2.3.

\f{\bf Proof of Lemma 2.3.}
The case of $k=0$ is induced from [Kf]. We prove the case of $k\geqq1$.

We first prove:

\f{\bf Claim.}
Let $F$ be a Seifert hypersurface for $K$. 
We can take $2p$ $(k+1)$-cycles  
$x_1,...,x_{p}$, $y_1,...,y_{p}$  
in $F$ which represent  basis of $H_{k+1}(F:\Bbb Z)$/(Torsion part) 
such that 
(1) $x_i\cdot y_i=1$ for all $i$,  
(Hence,$y_i\cdot x_i=(-1)^{k+1}$ for all $i$,)
(2)$x_i\cdot x_j=0$, $y_i\cdot y_j=0$ for all $(i,j)$, and 
(3)$x_i\cdot y_j=0$ for $i\neq j$.

Proof.
When $k$ is even,  take symplectic basis.
When $k$ is odd, we need the following.

\f
{\bf Sublemma.} {\sl(See e.g. [S].)
A symmetric matrix satisfies the conditions that 
(1)the elements are integers, 
(2)the  determinant is $\underline{+}1$, and 
(3)the  signature is zero, 
then the matrix is congruent to 
$\oplus \pmatrix
0&1\\
1&0\\
\endpmatrix.$}

\f Since a matrix which represents the intersection products 
of $H_{k+1}(F;\Bbb Z)$ satisfies the condition of the above sublemma,  
such $x_i$ and $y_j$ exist. 
The proof of the above Claim is completed.

There exists a Seifert matrix $X$ of $K$ 
associated with the basis, $x_i$ and $y_j$, and $F$. 
The elements of $X$ are $\theta(x_i,x_j)$ 
 $\theta(x_i,y_j)$,  
 $\theta(y_i,x_j)$, and
 $\theta(y_i,y_j)$ ($i,j=1,...,p$).

Take an embedding $f:S^{k+1}\x S^{k+1}\e B^{2k+3}$  
so that $f(S^{k+1}\x S^{k+1})$ is a tubular neighborhood of 
the standard ($k+1$)-sphere embedded trivially in $B^{2k+3}$.  
We regard $S^{k+1}$ as $D^{k+1}_1\cup D^{k+1}_2$. 
Let $A$ denote 
$f(S^{k+1}\x S^{k+1}-\roman{Int}\{D^{k+1}_2\x D^{k+1}_2\})$.  
Let $p_1,...,p_\mu$ be points in $D^{k+1}_2$, 
where $\mu$ is a large positive integer. 
Take a neighborhood 
$\cases
U_\alpha\\
V_\beta\\
\endcases$
of 
$\cases
D^{k+1}_1\x p_\alpha\\
p_\beta\x D^{k+1}_1\\
\endcases$
in $B^{2k+3}$ such that 
(1)$U_i$ and $V_i$ are diffeomorphic to open ($2k+3$)-balls, and 
(2)arbitrary two of them don't intersect. 
Let $q$ be the center of $D^{k+1}_1$.
Let $x'$ (resp.$y'$) denote the homology class which is represented by 
$q\x S^{k+1}$ (resp.$S^{k+1}\x q$). 
Note that 
(1)A cycle representing 
$x'$ intersects with each $U_\alpha$ at one points 
and doesn't intersect any $V_\beta$, 
and  
(2)A cycle representing 
$y'$ intersects with each $V_\beta$ at one points 
and doesn't intersect any $U_\alpha$.

Take disjoint ($2k+3$)-balls $B^{2k+3}_i(i=1,...,p)$ in $S^{2k+3}.$  
Take a copy of $A$ in each $B^{2k+3}_i$, say $A_i$. 
Take a copy of $U_\alpha$ (resp.$V_\beta$) in $B_i$, 
say $U_{i\alpha}$ (resp. $V_{i\beta}$).
Take a copy of $x'$ (resp.$y'$) in $B_i$, 
say $x'_i$ (resp. $y'_i$).
By using (2k+3)-dimensional 1-handles, 
take the connected sum of $A_i$ in $S^{2k+3}$, say $A_0$. 
Then $\partial A_0$ is the trivial ($2k+1$)-knot.

There exists a Seifert matrix $X'$ of the trivial knot  
associated with the basis, 
$x'_1,y'_1,...,x'_p,y'_p$,   
 and $A_0$. 
Obviously, $X'$ is 
$\oplus \pmatrix 
0&1\\
0&0\\
\endpmatrix.$
The elements of $X'$ satisfies that 
(1) $\theta(x'_i,x'_j)$=0 for all ($i,j$), 
(2) $\theta(y'_i,y'_j)$=0 for all ($i,j$), 
(3)  $\theta(x'_i,y'_j)$=$\theta(y'_j,x'_i)$=0 for $i\neq j$, 
(4) $\theta(x'_i,y'_i)$=1  for all $i$,  and 
(5)$\theta(y'_i,x'_i)$=0  for all $i$  ($i,j=1,...,p$).

We make a ($2k+1$)-knot $\widetilde{K}$ from the trivial ($2k+1$)-knot 
by high dimensional pass-moves  as in the following paragraphs. 
Before making $\widetilde{K}$, we prove: 
\f{\bf Claim.} 
If a Seifert matrix of $\widetilde{K}$ coincides with that of $K$,  
then
$\widetilde{K}$ is what we required, i.e., 
the proof of Lemma 2.3 is completed.

\f
Proof.
By the definition of the construction of $\widetilde{K}$, 
 $\widetilde{K}$ is pass-move equivalent to the trivial knot.  
By Theorem 2.5, $\widetilde{K}$ is cobordant to $K$.

We take the following pass-move-charts of $\partial A_0$, 
carry out the pass-moves, 
and modify the Seifert matrix $X'$ to $X$.   
For the new knots obtained by the following pass-moves, 
we can take a Seifert hypersurfaces diffeomorphic to $A$. 
We can call basis for the new Seifert matrixes $x_i$ and $y_j$, again.

Step 1. See $\theta (x_i,x_j)$ ($i>j$).
If $\theta (x_i,x_j)$=0, then 
 $\theta (x'_i,x'_j)$=$\theta (x_i,x_j)$. 
If $\theta (x_i,x_j)$=$\nu\neq$0, 
 make $\vert\nu\vert$ pass-move-charts $\widetilde{U_{ijk}}$ 
 ($k=1,...,$ $\vert\nu\vert$)
from $U_{ik}$ and  $U_{jk}$    
so that pass-moves in  
 $\widetilde{U_{ijk}}$ let $\theta (x'_i,x'_j)$=$\theta (x_i,x_j)$. 

Here, we have the following. We prove: 

\f{\bf Claim.}
Then $\theta (x'_i,x'_j)$=$\theta (x_i,x_j)$  ($i<j$).

Proof. Let $i<j$.
 By Lemma 2.4(1),  
 $\theta (x'_i,x'_j)$= $x'_i\cdot x'_j$+ $(-1)^{k}$ $\theta (x'_j, x'_i)$ and 
 $\theta (x_i,  x_j)$= $x_i\cdot x_j$+   $(-1)^{k}$ $\theta (x_j, x_i)$. 
 By the definition of  $x_i$ and $x'_j$,   
 $x_i\cdot x_j$=$x'_i\cdot x'_j$(=0).        
 Since $j>i$, 
 $\theta (x'_j,x'_i)$=$\theta (x_j,x_i)$. 
 Therefore 
 $\theta (x'_i,x'_j)$=$\theta (x_i,x_j)$.

Here, note that, by the definition of these pass-moves, 
each pass-move  in $\widetilde{U_{ijk}}$
doesn't change the value of $\theta (*,\dagger)$ except for 
$\theta (x'_i,x'_j)$ and $\theta (x'_j,x'_i)$ ($i>j$).

Step 2. See $\theta (y_i,y_j)$ ($i>j$).
If $\theta (y_i,y_j)$=0, then 
 $\theta (y'_i, y'_j)$=$\theta (y_i,y_j)$. 
If $\theta (y_i, y_j)$=$\nu\neq$0, 
 make $\vert\nu\vert$ pass-move-charts $\widetilde{V_{ijk}}$ 
 ($k=1,...,$ $\vert\nu\vert$)
from $V_{ik}$ and  $V_{jk}$    
so that pass-moves in  $\widetilde{V_{ijk}}$  
let $\theta(y'_i,y'_j)$=$\theta (y_i,y_j)$. 
Here, $\theta (y'_i,y'_j)$=$\theta (y_i,y_j)$  ($i<j$) holds by  Lemma 2.4(1).

Here, note that, by the definition of these pass-moves, 
each pass-move  in $\widetilde{V_{ijk}}$
doesn't change the value of $\theta (*,\dagger)$ except for 
$\theta (y'_i, y'_j)$ and $\theta (y'_j, y'_i)$ ($i>j$).

Step 3. See $\theta (x_i,y_j)$ for any $(i,j)$.
If $\theta (x_i,y_j)$=0, then 
 $\theta (x'_i, y'_j)$=$\theta (x_i,y_j)$. 
If $\theta (x_i, y_j)$=$\nu\neq$0, 
 make $\vert\nu\vert$ pass-move-charts $\widetilde{W_{ijk}}$ 
 ($k=1,...,$ $\vert\nu\vert$)
from $U_{ik}$ and  $V_{jk}$    
so that pass-moves in  $\widetilde{W_{ijk}}$  
let $\theta(x'_i, y'_j)$=$\theta (x_i, y_j)$. 
Here,  $\theta (y'_i, x'_j)$=$\theta (y_i, x_j)$   holds  by  Lemma 2.4(1).

Here, note that, by the definition of these pass-moves, 
each pass-move in  $\widetilde{W_{ijk}}$
doesn't change the value of $\theta (*,\dagger)$ except for 
$\theta (x'_i, y'_j)$ and $\theta (y'_j, x'_i)$.


Here, note that, by the definition of these pass-moves, 
each pass-move  in $\widetilde{W'_{ijk}}$  
  doesn't change the value of $\theta (*,\dagger)$ except for 
$\theta (y'_i, x'_j)$ and $\theta (y'_j, x'_i)$ ($i>j$).

Before Step 4 and 5, we prove:

\f{\bf Claim.}
(1)If $k$ is odd,
 $\theta (x_i,x_i)$=$\theta (y_i,y_i)$=0.   
(2)If $k$ is even,we can assume 
 $\theta (x_i,x_i)$, and $\theta (y_i,y_i)$ 
 are even integer.

  Proof. 
  (1)By  Lemma 2.4(1) and the definitions of 
  $x_i$, $y_j$, $x'_i$ and $y'_j$, 
 $\theta (x'_i,x'_i)$=$\frac{x'_i\cdot x'_i}{2}$=0 and 
 $\theta (y'_i,y'_i)$=$\frac{y'_i\cdot y'_i}{2}$=0.  
   (2) We can change the basis, if necessary, because Arf$(K)$=0 
   and the Arf invariant of the trivial knot is zero.

Step 4. 
See $\theta (x_i, x_i)$.   
If $k$ is odd,  $\theta (x_i,x_i)$=0.   
Hence  $\theta (x'_i, x'_i)$=$\theta (x_i,x_i)$. 
 The case when $k$ is even.  
If $\theta (x'_i,x'_i)$=0, then 
 $\theta (x'_i, x'_i)$=$\theta (x_i,x_i)$. 
The case of  $\theta (x'_i, x'_i)$$\neq$0. 
We can put it $2\nu$. 
 Make $\vert\nu\vert$ pass-move-charts $\widetilde{Z_{ik}}$ 
 ($k=1,...,$ $\vert\nu\vert$)
from $U_{ik}$ and  $U_{i(\nu+k)}$    
so that pass-moves in  $\widetilde{Z_{ik}}$  
let $\theta(x'_i,x'_i)$=$\theta (x_i,x_i)$. 

Here, note that, by the definition of these pass-moves, 
each pass-move in  $\widetilde{Z_{ik}}$ 
doesn't change the value of $\theta (*,\dagger)$ except for 
$\theta (x'_i,x'_i)$.

Step 5.(final step.) 
See $\theta (y_i,y_i)$.   
If $k$ is odd,  $\theta (y_i,y_i)$=0. 
Hence $\theta (y'_i, y'_i)$=$\theta (y_i,y_i)$. 
 The case $k$ is even.  
If $\theta (y_i,y_i)$=0, then 
 $\theta (y'_i, y'_i)$=$\theta (y_i,y_i)$. 
The case of  $\theta (y_i, y_i)$$\neq$0. 
We can put it $2\nu$. 
 Make $\vert\nu\vert$ pass-move-charts $\widetilde{Z'_{ik}}$ 
 ($k=1,...,$ $\vert\nu\vert$)
from $V_{ik}$ and  $V_{i(\nu+k)}$    
so that pass-moves in  $\widetilde{Z'_{ik}}$  
let $\theta(y'_i, y'_i)$=$\theta (y_i, y_i)$.

Here, note that, by the definition of these pass-moves, 
each pass-move in $\widetilde{Z'_{ik}}$  
doesn't change the value of $\theta (*,\dagger)$ except for 
$\theta (y'_i,y'_i)$.

 We now obtain $\widetilde{K}$. 
As seen before, we complete the proof of Lemma 2.3. 
\qed


Note that, in the above proof, we proved the following Claim 2.2.1.S 
which is stronger than Claim 2.2.1.

\f{\bf Claim 2.2.1.S.} {\sl 
If (2) of Theorem 2.2 holds, then (1) of Theorem 2.2 holds and 
there exist pass-move-charts $U_i$ $(i=1,..,q)$ of $K_3$ such  that 
$U_i$$\cap$$U_j$=$\phi$ $(i\neq j)$ and 
$K_1$ is obtained from $K_3$ by the pass-moves in $U_i$.  }



\head 3.
A sufficient condition for the realization of pair of odd dimensional knots  
\endhead

In this section we prove: 
\proclaim{ Proposition 3.1}
Let $K_1$ and $K_2$ be ($2k+1$)-knots and 
$X_2$  a ($2k+3$)-knot   
diffeomorphic to the standard ($2k+3$)-sphere ($k\geqq0$). 
 If we have the condition that 
 \nl
 $\cases
\text{ 
 Arf$(K_1)$=Arf$(K_2)$ 
 }
 & \text{when $k$ is even}\\ 
 \text{ 
$\sigma (K_1)=\sigma (K_2)$ 
 }
 & \text{when $k$ is odd,} 
\endcases$
then, for a ($2k+3$)-knot  $X_1$ diffeomorphic to the standard 
($2k+3$)-sphere,  ($K_1$, $K_2$,$X_1$, $X_2$) is realizable. 
\endproclaim

To prove  Proposition 3.1, we need some lemmas. 

\proclaim{Lemma 3.2}
Let  $K$ be an arbitrary $n$-knot and  
$X_1$ and $X_2$  arbitrary $(n+2)$-knots diffeomorphic to 
the standard $(n+2)$-sphere ($n\geqq1$). 
Then  $(K,K,X_1, X_2)$ is realizable.
\endproclaim

\proclaim{Lemma 3.3}
If 4-tuple of $(n,n+2)$-knots 
($K_1$, $K_2$, $X_1$, $X_2$) is realizable ($n\geqq1$) and  
 $K_2$ is cobordant to $\widetilde{K_2}$, 
 then for an $(n+2)$-knot $\widetilde{X_1}$  diffeomorphic to  
 the standard $(n+2)$-sphere, 
 ($K_1$, $\widetilde{K_2}$, $\widetilde{X_1}$,  $X_2$) is realizable.  
 Furthermore, for an arbitrary Seifert hypersurface $F$ 
 for the  $n$-knots $\widetilde{K_2}$, 
 there exists an immersion 
$\widetilde{f}:S^{n+2}_1\coprod S^{n+2}_2 \looparrowright S^{n+4}$
 to realize 
 ($K_1$, $\widetilde{K_2}$, $\widetilde{X_1}$,  $X_2$)  
 and a Seifert hypersurface $\widetilde{V}$ for 
 $\widetilde{X_1}$=$\widetilde{f}(S^{n+2}_1)$ such that 
 $\widetilde{V}$ $\cap$ $\widetilde{f}(S^{n+2}_2)$= 
 $\widetilde{V}$ $\cap$ $X_2$ 
 is the Seifert hypersurface $F$ for $\widetilde{K_2}$.   
\endproclaim

\proclaim{ Lemma 3.4}
Let $K_1$and $K_2$ be $(2k+1)$-knots and 
$X_1$ and $X_2$  arbitrary $(2k+3)$-knots diffeomorphic to 
the standard $(2k+3)$-sphere ($k\geqq0$). 
If $K_1$and $K_2$ are pass-move equivalent and 
there exist pass-move-charts $U_i$ $(i=1,..,q)$ of $K_2$ such  that 
$U_i$$\cap$$U_j$=$\phi$ $(i\neq j)$ and 
$K_1$ is obtained from $K_2$ by the pass-moves in $U_i$,   
then ($K_1$, $K_2$, $X_1$, $X_2$) is realizable.  
\endproclaim

We first prove: 

\f{\bf Claim.}
If 
Claim 2.2.1.S, 
Lemma 3.3, and 3.4  hold,  
Proposition 3.1 holds.

Proof. 
By  Claim 2.2.1.S, 
(1) there exists a $(2k+1)$-knot $K_3$  
which is pass-move equivalent to $K_1$ and cobordant to $K_2$ and 
(2)there exist pass-move-charts $U_i$ $(i=1,..,q)$ of $K_3$ such  that 
$U_i$$\cap$$U_j$=$\phi$ $(i\neq j)$ and 
$K_1$ is obtained from $K_3$ by the pass-moves in $U_i$.   
By Lemma 3.4, ($K_1, K_3, X_2, X_2$) is realizable. 
By Lemma 3.3, for a ($2k+3$)-knot $X_1$, 
($K_1, K_2, X_1, X_2$) is realizable.

Note.  Lemma 3.2 is used to prove Lemma 3.3 and 3.4.
The latter half of Lemma 3.3 is used in \S 4. 
The case when $n$ is even of Theorem 1.3 is induced from Lemma 3.4, obviously. 
But we prove in \S 6 Theorem 1.4 stronger than it.

In the rest of this section we prove Lemma 3.2-4.

\f{\bf Proof of Lemma 3.2. }
Take an embedding $f_a:S^{n+2}_1\coprod S^{n+2}_2 \e S^{n+4}$ 
which defines the trivial $(n+2)$-link.
There exists a chart $U$ of $S^{n+4}$ 
with the following properties (1) and (2). 
\roster
\item
 $\phi:U\cong \Bbb R^{n+2}$ $\times$
$\{(u,v)\vert u,v\in  \Bbb R\}$ 
$\cong \Bbb R^{n+2}$ $\times \Bbb R_u \times \Bbb R_v$

\item
$U\cap f_a(S^{n+2}_1)
=\Bbb R^{n+2}\times\{(u,v)\vert u=0, v=0\}$

\f $U\cap f_a(S^{n+2}_2)
=\Bbb R^{n+2}\times\{(u,v)\vert u=3, v=0\}$
\endroster

We modify the embedding $f_a$ to 
obtain an immersion  
$f_b:S^{n+2}_1\coprod S^{n+2}_2 \looparrowright S^{n+4}$ 
 to realize $(K, K, T, T)$, where $T$ is the trivial knot.   
 We put $f_b\vert S^{n+2}_2=f_a\vert S^{n+2}_2$. 
 We define  $f_b\vert S^{n+2}_1$ as follows. 
 Take an embedding $g:\Sigma^{2k+1}\e$ $U\cap f_a(S^{n+2}_1)$ 
which defines the  $n$-knot $K$ in $S^{n+2}_1$.
Let $F$ be a Seifert hypersurface for $K$ and 
a submanifold of $U\cap f_a(S^{n+2}_1)$. 
Let   $N_1(F)=F\times[-1,1]$ be a submanifold embedded  
in  $U\cap f_a(S^{n+2}_1)$ such that $F=F\times\{0\}$. 
We define the subsets  $E_1, E_2 $ and $E_3$ of 
$N_1(F) \times \Bbb R_u \times \Bbb R_v =
\{(p,t,u,v)\vert p\in F, -1\leqq t\leqq 1 , u\in \Bbb R , v\in \Bbb R \}$
 as follows. 

\hskip3mm$E_1=\{(p,t,u,v)\vert p\in F, 0\leqq u\leqq 1, -1\leqq t\leqq 1, v=0\}$

\hskip3mm$E_2=\{(p,t,u,v)\vert p\in F, 1\leqq u\leqq 2, 
t=k\cdot \roman{cos}{\frac {\pi (u-1)}{2}}, 
v=k\cdot \roman{sin}{\frac {\pi (u-1)}{2}}, 
-1\leqq k\leqq 1\}$ 

\hskip3mm$E_3=\{(p,t,u,v)\vert p\in F,2\leqq u\leqq 4,t=0,$ 
$-1\leqq v\leqq 1 \}$

Then the followings hold by the way of the construction. 

 $\Sigma=$
$\overline{\text
{$\{f_a(S^{n+2}_1)-N_1(F)\}$}}$  
 $\cup_{ \partial N_1(F)}$
$\overline{\text
\{\partial (E_1\cup E_2\cup E_3)\}-N_1(F)}$
\nl is a (n+2)-sphere embedded in $S^{n+4}$. 
\nl Since 
$\overline{\text
\{\partial (E_1\cup E_2\cup E_3)\}-N_1(F)}$
is isotopic to $N_1(F)$ relative to $\partial N_1(F)$, 
$f_b\vert S^{n+2}_1$ defines the trivial knot. 
 Both  
 $\Sigma\cap$ $f_b(S^{n+2}_2)$  in  $\Sigma$ and 
 $\Sigma\cap$ $f_b(S^{n+2}_2)$  in  $f_b(S^{n+2}_2)$  
 defines $K$. 
We define $f_b\vert S^{n+2}_1$ so that $f_b(S^{n+2}_1)$ 
coincides with $\Sigma$. 
We obtain $f_b$  to realize $(K, K, T, T)$.


By the following Sublemma 3.5, $(K, K,$$ X_1,X_2)$ is realizable. 
We prove Sublemma 3.5 and complete the proof of Lemma 3.2.

{\bf Sublemma 3.5. }
Let $T$ be the trivial ($n+2$)-knot and  $X_i$ arbitrary  ($n+2$)-knots 
($i=1,2$)($n\geqq1$). 
If $(K_1, K_2, T, T)$ is realizable, 
then  $(K_1, K_2,X_1, X_2)$ is realizable. 

{\p Sublemma 3.5.}Let $f':S^{n+2}_1\amalg S^{n+2}_2 \i S^{n+4}$ realize $(K_1, K_2, T, T)$. Let $f_i:S^{n+2}_i \e S^{n+4} (i=1,2)$ define $(n+2)$-knots $X_i$. Let $B', B_1$ and $B_2$ be ($n+4$)-balls in $S^{n+4}$ and $B'\cap B_1=B'\cap B_2=B_1\cap B_2=\phi$.  We can take $f'$ and $f_i$  so that Im$f_i$ in $B_i$ and Im$f'$ in $B'$. Connect  $f_i(S^{n+2}_i)$ with   $f'(S^{n+2}_i)$ by ($n+3$)-dimensional 1-handle $h_i$ embedded in $S^{n+4}$ ($i=1,2$), where $h_1\cap h_2=\phi.$Take $f:S^{n+2}_1\amalg S^{n+2}_2 \i S^{n+4}$ so that  $f(S^{n+2}_i)$ coincides with   $f_i(S^{n+2}_i)$ $\sharp$  $f'(S^{n+2}_i)$. Then $f$ realizes $(K_1, K_2,X_1, X_2)$. \qed

\f{\bf Proof of Lemma 3.3. }
Let ${f}:S^{n+2}_1\coprod S^{n+2}_2 \looparrowright S^{n+4}$
be  an immersion  to realize $(K_1, K_2, X_1, X_2)$ 
and  $V$ a Seifert hypersurface  for ${f}(S^{n+2}_1)$. 
Let ${f}(S^{n+2}_2)$$\times$ $D^2$=$X_2\times$
$\{(x,y) \vert$ $x=r\cdot\co\theta, y=r\cdot\si\theta,$  
$0\leqq r\leqq1, 0\leqq \theta <2\pi  \}$
be a tubular neighborhood of $X_2$  in $S^{n+4}$. 

See 
 ${V}$ $\cap$ $\{f(S^{n+2}_2)$ $\times$ $D^2\}$. 
It has the following properties. 
\roster
\item
For any $(x, y)$,  
 $\{\partial\widetilde{V}\}$ $\cap$ $[ f(S^{n+2}_2)$ $\times$ $(x, y)$]
 in $f(S^{n+2}_2)$ $\times$ $(x, y)$  
 defines  ${K_2}$.
\item
For each $(x, y)$, 
$\widetilde{V}$  $\cap$ $[ f(S^{n+2}_2)$ $\times$ $(x,y)$] 
is a  same Seifert hypersurface $G$.  
 
\item
For any $\theta$, 
 $\{\partial{V}\}$ $\cap$ 
 $[ f(S^{n+2}_2)$ $\times$  
 $\{(x,y)\vert x=r\cdot\co\theta, y=r\cdot\si\theta, 0\leqq r \leqq1\} $]
 is diffeomorphic to $K_2\times[0,1]$. 

\endroster

Prepare $S^{n+2}$ $\times$ $[0,1]$. 
Put $K_2$ and $G$ in $S^{n+2}$ $\times$ $\{1\}$. 
Put $\widetilde{K_2}$ and $F$ in $S^{n+2}$ $\times$ $\{0\}$. 
Recall $K_2$ and  $\widetilde{K_2}$ are cobordant. 
Hence there exists a submanifold $P$ in  $S^{n+2}$ $\times$ $[0,1]$  
which meets the boundary transversely in $P$, 
is diffeomorphic to $K_2$$\x$[0,1], 
and meets  
$S^{n+2}$ $\times$ $\{0\}$  at $K_2$ and 
$S^{n+2}$ $\times$ $\{1\}$  at  $\widetilde{K_2}$. 
By an elementary discussion of the obstruction theory, 
there exists a compact submanifold $Q$ in 
 $S^{n+2}$ $\times$ $[0,1]$ such that 
 $\partial Q$=$F$ $\cup$ $P$ $\cup$ $G$.

We modify $V$ and make the following  $\widetilde{V}$.  
Let $\widetilde{V}$ be a submanifold of $S^{n+4}$ satisfying 
the following conditions. 
\roster

\item
 $\widetilde{V}$ $\cap$ $[S^{n+4}-$ $\{f(S^{n+2}_2)$ $\times$ $D^2\}$] 
 coincides with 
 ${V}$ $\cap$ $[S^{n+4}$- $\{f(S^{n+2}_2)$ $\times$ $D^2\}]$. 

\item
 $\{\partial\widetilde{V}\}$ $\cap$ $[ f(S^{n+2}_2)$ $\times$ $(0,0)$]
 in $f(S^{n+2}_2)$ $\times$ $(0,0)$
 is $\widetilde{K_2}$.

\item
$\widetilde{V}$  $\cap$ $[ f(S^{n+2}_2)$ $\times$ $(0,0)$] =$F$. 

\item 
For any $\theta$, 
 $\{\partial\widetilde{V}\}$ $\cap$ 
 $[ f(S^{n+2}_2)$ $\times$  $(\co\theta, \si\theta)$]
  in  
 $[ f(S^{n+2}_2)$ $\times$  $(\co\theta, \si\theta)$ ]
  is $K_2$.
 
\item
For any $\theta$, 
 $\{\partial\widetilde{V}\}$ $\cap$ 
 $[ f(S^{n+2}_2)$ $\times$  
 $\{(x,y)\vert x=r\cdot\co\theta, y=r\cdot\si\theta, 0\leqq r \leqq1\} $]
 is diffeomorphic to $K_2\times[0,1]$( and $\widetilde{K_2}$$\times[0,1]$). 
\item 
For any $\theta$, 
 $\{\widetilde{V}\}$ $\cap$ 
 $[ f(S^{n+2}_2)$ $\times$  
 $\{(x,y)\vert x=r\cdot\co\theta, y=r\cdot\si\theta, 0\leqq r \leqq1\} $]
 defines the above submanifold $Q$ in 
 $[ f(S^{n+2}_2)$ $\times$  
 $\{(x,y)\vert x=r\cdot\co\theta, y=r\cdot\si\theta, 0\leqq r \leqq1\} $].

\endroster
Note.(i) The above condition (5) holds because 
   $K_2$ is cobordant to $\widetilde{K_2}$.
(ii)$\partial\widetilde{V}$ is diffeomorphic to the standard sphere.    

Let $\widetilde{f}:S^{n+2}_1\coprod S^{n+2}_2 \looparrowright S^{n+4}$
be an immersion such that 
(1) $\widetilde{f}(S^{n+2}_1)$ coincides with $\partial\widetilde{V}$, 
 say $\widetilde{X_1}$, and 
(2) $\widetilde{f}\vert S^{n+2}_2={f}\vert S^{n+2}_2$.

By the construction of $\widetilde{f}$, 
the $n$-knot 
 $\widetilde{f}(S^{n+2}_1)$ $\cap$  $\widetilde{f}(S^{n+2}_2)$ 
 in  $\widetilde{f}(S^{n+2}_1)$ 
is equivalent to 
 $\widetilde{f}(S^{n+2}_1)$ $\cap$  
 $[{f}(S^{n+2}_2)\x \{(x,y)\vert x=1,y=0\}]$ 
 in  $\widetilde{f}(S^{n+2}_1)$,  
\nl to 
 ${f}(S^{n+2}_1)$ $\cap$  
 $[{f}(S^{n+2}_2)\x \{(x,y)\vert x=1,y=0\}]$ 
 in  ${f}(S^{n+2}_1)$,  
\nl and to 
 ${f}(S^{n+2}_1)$ $\cap$  ${f}(S^{n+2}_2)$  in  ${f}(S^{n+2}_1)$, 
 that is, $K_1$.

\f
Then we obtain 
the immersion $\widetilde{f}$ to realize  
 ($K_1$, $\widetilde{K_2}$, $\widetilde{X_1}$,  $X_2$)  
and  the required $V$.




\f{\bf Proof of Lemma 3.4. }
 Take the pass-move-charts $U_1,...,U_q$ of $K_2$. 
 In each pass-move-chart $U_i$, take 

\hskip2cm
 $D^{k+1}_{1i}$= 
 $[\frac{7}{18},\frac{11}{18}] 
 \times D^{k+1}_x(0)
 \times D^{k+1}_y(\frac{1}{2})$  

\hskip2cm
 $D^{k+1}_{2i}$=
 $[\frac{5}{9},\frac{7}{9}   ] 
 \times D^{k+1}_x(\frac{1}{2}) 
 \times D^{k+1}_y(0)$. 

\f Put  $S^{k+1}_{1i}$= $\partial D^{k+1}_{1i}$
and  $S^{k+1}_{2i}$= $\partial D^{k+1}_{2i}$. 
Note the linking number of $S^{k+1}_{1i}$ and $S^{k+1}_{2i}$  is one. 
Since $H_{2k+1}(S^{2k+3}-\cup_{i,j}\{S^{k+1}_{ji}\})$=0, 
a Seifert hypersurface $F$ for $K_2$ is included in 
$S^{2k+3}-\cup_{i,j}S^{k+1}_{ji}$ ($i=1,...,q, j=1,2 $).

Then  we have the following.

\f{\bf Claim.}   {\sl  
If we attach ($2k+4$)-dimensional ($k+2$)-handles with 0-framing 
to $U_i$ along 
$S^{k+1}_{1i}$   and $S^{k+1}_{2i}$
and carry out surgery on $U_i$ (and $S^{2k+3}$), 
then 
(1)$U_i$       changes into the $(2k+3)$-ball again, 
(2)$S^{2k+3}$  changes into the $(2k+3)$-sphere  again, and 
(3)the new knot in the new $(2k+3)$-sphere is 
the knot obtained from $K_2$ by the pass-moves in $U_i$. 
}


By the above discussion, the followings hold. 
In all $U_i$ ($i=1,...,q$), carry out the above surgeries. 
Then 
$S^{2k+3}$  changes into the $(2k+3)$-sphere  again, and 
$K_2$ in the old $(2k+3)$-sphere  changes into 
$K_1$ in the new $(2k+3)$-sphere. 
There exists a Seifert hypersurface $F$ for $K_2$ 
such that 
$S^{2k+3}$ $\cap$ $S^{k+1}_{ji}$ for all $i,j$.


We first construct an immersion to realize 
 $(K_2,K_2,T,T)$ as in the proof of Lemma 3.2. 
Take $U$, and $\Sigma$  as in the proof of Lemma 3.2.  
 Take the pass-move-charts $U_i$ of $K_2$ in $U$. 
Take  $S^{k+1}_{ji}$ in $U_i$ $\subset$ $U$.  
We use the Seifert hypersurface $F$ in the previous paragraph 
as $F$ in the proof of Lemma 3.2. 
 As we see before,  we can take a Seifert hypersurface  
 so that it does not intersect with  $S^{k+1}_{ji}$ in $U_i$ for all $i, j$.
We use the Seifert hypersurface as $F$  in the proof of Lemma 3.2.  
We use these $U$,  $\Sigma$ and $F$ and 
construct an immersion 
 $f_b:S^{2k+3}_1\coprod S^{2k+3}_2 \i S^{2k+5}$   
 to realize   $(K_2,K_2,T,T)$ as in the proof of Lemma 3.2.

We next construct an immersion 
$f_c:S^{2k+3}_1\coprod S^{2k+3}_2 \e S^{2k+5}$
to realize ($K_1, K_2$,T,T). 
 Let $h_{1i}$ be a tubular neighborhood of 
 
[$S^{k+1}_{1i}\x \{(u,v)\vert u=0, 0\leqq v \leqq 1\}$]
 $\cup$ 
[$D^{k+1}_{1i}\x \{(u,v)\vert u=0, v=1\}$] 
\nl in $U_i\x \{(u,v)\vert u=0,  v\geqq 0 \}$. 
Let $h_{2i}$ be a tubular neighborhood of 

$[S^{k+1}_{2i}\x \{(u,v)\vert u=0, -1\leqq v \leqq 0\}]$ $\cup$
$[D^{k+1}_{2i}\x \{(u,v)\vert u=0, v=-1\}]$ 
\nl in $U_i\x \{(u,v)\vert u=0,  v\leqq 0 \}$. 
Make a submanifold $\Lambda$ from $\Sigma$ and $h_{ji}$ so that 

 $\Lambda$=
 $\overline{
 \text{
 $\Sigma-$
 $\cup^q_{i=1}$
 $(h_{1i}\cap\Sigma)-$
 $\cup^q_{i=1}$
 $(h_{2i}\cap\Sigma)$
 }
 }$
 $\cup^q_{i=1}$
 $\overline{
 \text{
 $\partial h_{1i}- (h_{1i}\cap\Sigma$)  
 }
 }$
 $\cup^q_{i=1}$
 $\overline{
 \text{
 $\partial h_{2i}- (h_{1i}\cap\Sigma$) 
 }
 }$. 
\nl  Of course,  $(h_{ji}\cap\Sigma)$ is  $(\{\partial h_{ji}\}\cap\Sigma)$ 
 and 
 is the tubular neighborhood of  $S^{k+1}_{ji}$ in  $U_i$.


Then  the followings hold by the definition  of the construction.
(1)$\Lambda$ is the trivial $(2k+3)$-knot.  
(2)$\Lambda\cap$ $f_b(S^{2k+3}_2)$  in  $\Lambda$ is $K_1$. 
(Because the pass-move is carried out in each pass-move-chart $U_i$.) 
(3)$\Lambda\cap$ $f_b(S^{2k+3}_2)$   in $f_b(S^{2k+3}_2)$ is  $K_2$.  
Here, we take an immersion 
$f_c:S^{2k+3}_1\coprod S^{2k+3}_2 \e S^{2k+5}$
so that 
$f_c(S^{2k+3}_1)$ coincides with $\Lambda$ and 
$f_c\vert S^{2k+3}_2 =f_b\vert S^{2k+3}_2$. 
Then $f_c$ is an immersion to realize ($K_1, K_2$,$T,T)$. 

At last, by Sublemma 3.5,  ($K_1, K_2$,$X_1, X_2$) is realizable. 
\qed


\head 4.
A necessary condition 
for the realization of pair of $(4m+1)$-knots  
\endhead


In this section we prove: 
\proclaim{Proposition 4.1}
If a 4-tuple of ($4m+1, 4m+3$)-knots ($K_1, K_2, X_1, X_2$) is realizable 
($m\geqq0$), then Arf($K_1$)=Arf($K_2$).
\endproclaim

We first prove the following 
Lemma 4.2.
\f
{\bf Lemma 4.2}  {\sl 
If 4-tuple of $(n, n+2)$-knots 
($K^{(1)}_1,$ $ K^{(1)}_2,$ $ X^{(1)}_1,$ $ X^{(1)}_2$) 
and 
($K^{(2)}_1,$ $ K^{(2)}_2,$ $ X^{(2)}_1,$ $ X^{(2)}_2$) 
are realizable ($n\geqq1$), 
then 
($K^{(1)}_1\sharp $ $K^{(2)}_1,$$ K^{(1)}_2\sharp$ $ K^{(2)}_2$, 
 $X^{(1)}_1\sharp $ $X^{(2)}_1,$$ X^{(1)}_2\sharp$ $ X^{(2)}_2$) 
is realizable. }

Note. Lemma 4.2 includes Sublemma 3.5.

{\p  Lemma 4.2.}
Take two  $(n+4)$- balls $B^{n+4}_1$ and $B^{n+4}_2$ in $S^{n+4}$ so that 
$B^{n+4}_1\cap B^{n+4}_2$=$\phi$. 
Take $g_i:S^{n+2}_{1}\amalg S^{n+2}_{2} \i S^{n+4}$ to 
realize 
($K^{(i)}_1, K^{(i)}_2, X^{(i)}_1, X^{(i)}_2$) 
so that we take Im$g_i$  in $B^{n+4}_i (i=1,2)$. 
Let $X_{j}^{(i)}$ denote Im$g_{i}(S^{n+2}_{j})$ for conveinience.
Connect $X^{(1)}_1$ with $X^{(2)}_1$  
by using ($n+3$)-dimensional 1-handle $h_{g1}$ embedded in $S^{n+4}$ 
to obtain  $\widetilde{X_1}$= $X^{(1)}_1\sharp X^{(2)}_1$. 
Connect $\widetilde{X_1}$$\cap$ $X^{(1)}_2$ 
with $\widetilde{X_1}$$\cap$$X^{(2)}_2$ 
by ($n+1$)-dimensional 1-handle $h'_{g2}$ embedded in  $\widetilde{X_1}$. 
Connect $X^{(1)}_2$ with $X^{(2)}_2$ 
by using ($n+3$)-dimensional 1-handle $h_{g2}$= $h'_{g2}\x D^2$ embedded 
 in $S^{n+4}$ to obtain  $\widetilde{X_2}$= $X^{(2)}_2\sharp X^{(2)}_2$. 
Take $g:S^{n+2}_1\amalg S^{n+2}_2 \i S^{n+4}$ 
so that  $g(S^{n+2}_i)$ coincides with $\widetilde{X_i}$. 
\qed

We prove Proposition 4.1 by the reduction to absurdity. 
We assume Arf($K_1$)$\neq$ Arf($K_2$) 
and induce the absurdity.
We may assume that Arf($K_1$)=1 and  Arf($K_2$)=0 
without loss of generality.

Note. 
If $bP_{4m+2}$=$\Bbb Z_2$ , it is obvious that the proof is easy. 
Recall that 
the subgroup $bP_{4m+2}$ $\subset$ $\Theta_{4m+1}$ 
is the trivial group for some integers 
and 
is  $\Bbb Z_2$ for some integers (See [KM].) 
It is known that 
there exist some integers  $m$ 
such that 
 (i) $TS^{2m+1}$ is not the trivial bundle, i.e., $m\neq$0,1,3, and 
 (ii) $bP_{4m+2}$ is the trivial group (See [Br]).

Let $T$ be the trivial $(4m+1)$-knot.  
Let $K$ be a $(4m+1)$-knot such that   
a Seifert hypersurface  for $K$ is the plumbing $F$ of two copies of 
the tangent bundle of the standard $(2m+1)$-sphere   
and a Seifert matrix associated with $F$ is 
$\pmatrix
1&1\\
0&1\\
\endpmatrix.$
(See e.g. P. 162 of [LM] for the plumbing.)
Then Arf($K$)=1 and for an arbitrary non-vanishing ($2m+1$)-cycle $x$ of $F$, 
$\theta(x,x)$ is odd. 
(See \S 2 of this paper for the definition of $\theta(\ ,\ )$.)

We prove:

\f{\bf Claim.}  {\sl
The pair of ($4m+1$)-knots ($T,K$) is realizable
under the hypothesis of the reduction to absurdity.}

{\bf Proof.} 
By Proposition 3.1, ($-K^*_1, K$) and ($T, -K^*_2$) are 
realizable. 
By the hypothesis of the reduction to absurdity   
($K_1, K_2$) is realizable. 
By Lemma 4.2, 
($K_1\sharp$$\{-K^*_1\}\sharp T,$  $ K_2\sharp K $$\sharp\{-K^*_2\}$)
is realizable, i.e., 
($K_1\sharp \{-K^*_1\},$   $K_2\sharp K \sharp \{ -K^*_2\}$)
 is realizable. 
 Since $K_1\sharp \{-K^*_1\}$ is cobordant to the trivial knot 
 and $K_2\sharp K $  $\sharp \{ -K^*_2\}$  
  is cobordant to $K$, 
 ($T,K$) is realizable  by Lemma 3.3. 
 \qed

Let  $f:S^{4m+3}_1\coprod S^{4m+3}_2 \looparrowright S^{4m+5}$ 
be  an  immersion to realize $(T, K)$. 
By Lemma 3.3,
there exist Seifert hypersurfaces $V_i$ for  $f(S^{4m+3}_i)$ ($i=1,2$)
such that   $f(S^{4m+3}_1)\cap V_2$ is the ($4m+2$)-disk $D$ and 
 $V_1\cap f(S^{4m+3}_2)$ is the Seifert hypersurface $F$. 
Let $W$ denote $V_1\cap V_2$. Then $\partial W$ =$F\cup D.$
Then the following holds.We prove:

\f{\bf Claim.}   {\sl 
There exists a non-vanishing  $(2m+1)$-$\Bbb Z_2$-cycle $x$ 
in $F$ which is zero cycle in $W$, i.e., 
$x$ bounds a ($2m+2$)$-\Bbb Z_2$-chain $y$ in $W$.}

  \def\ri{\rightarrow}

{\bf Proof. }  
The natural inclusion 
$F\e$ $\partial W$ 
induces 
$ H_i(F; $$\Bbb Z_2)$  $\cong$ $ H_i(\partial W; $$\Bbb Z_2)$ ($i\neq4m+2$). 
Hence it suffices to prove that 
Ker
$\{ H_{2m+1}$$(\partial W; $$\Bbb Z_2)$ $\ri$ $H_{2m+1}$$( W; $$\Bbb Z_2)\}$ 
is not the trivial group.
Consider the exact sequence: 
$H_{i}(\partial W;\Bbb Z_2)$ 
 $\rightarrow$
$H_{i}(W ;\Bbb Z_2)$ 
 $\rightarrow$
$H_{i}$ $(W,$$ \partial W;\Bbb Z_2)$. 
We use the following  part ($*$) of the exact sequence:
$H_{2m+2}(\partial W;\Bbb Z_2)$ $\rightarrow$ 
$H_{2m+2}(W ;\Bbb Z_2)$  $\rightarrow$ 
$H_{2m+2}$ $(W,$$ \partial W;\Bbb Z_2)$ $\rightarrow$ 
$H_{2m+1}(\partial W;\Bbb Z_2)$  $\rightarrow$
$H_{2m+1}(W ;\Bbb Z_2)$  $\rightarrow$
$H_{2m+1}$ $(W,$$ \partial W;\Bbb Z_2)$  $\rightarrow$
$H_{2m}(\partial W;\Bbb Z_2)$.   
We have 
$H_{2m+2}$ $(\partial W; $$\Bbb Z_2)$ $\cong$ 
$H_{2m}$ $(\partial W; $$\Bbb Z_2)$ $\cong$0 and 
$H_{2m+1}$$(\partial W;$$\Bbb Z_2)$ 
$\cong$ $\Bbb Z_2$ $\oplus$$\Bbb Z_2$. 
 By Poincar\'e duality and the universal coeficient theorem,  
 the followings hold. 
(1)The $\Bbb Z_2$-rank of  $H_{2m+2}$ $(W; $$\Bbb Z_2$) 
and 
that of $H_{2m+1}$ $(W, \partial W; $$\Bbb Z_2)$ 
are same, put it $r$. 
(2) The $\Bbb Z_2$-rank of  $H_{2m+1}$$( W; $$\Bbb Z_2$) 
and 
that of $H_{2m+2}$$(W, $$\partial W; $$\Bbb Z_2)$ 
are same, put it $s$. 
Therefore  the sequence $(*)$ become as follows:
 0 $\ri$
 $\oplus^r $$\Bbb Z_2$ $\ri$
 $\oplus^s $$\Bbb Z_2$ $\ri$ 
 $\Bbb Z_2 $$\oplus\Bbb Z_2$ $\ri$
 $\oplus^s $$\Bbb Z_2$ $\ri$ 
 $\oplus^r $ $\Bbb Z_2$$\ri$ 0. 
 Therefore the $\Bbb Z_2$-rank of  Ker$\{ H_{2m+1}$$(\partial W; $$\Bbb Z_2)$ $\ri$
 $ H_{2m+1}$$( W; $$\Bbb Z_2)\}$ =$s-r$=1.
  \qed

The ($2m+1$)$-\Bbb Z_2$-cycle $x$ bounds a ($2m+2$)$-\Bbb Z_2$-chain $z$ 
in $S^{4m+3}_2$. 
Let $w$ denote the $(2m+2)$-$\Bbb Z_2$-cycle $y\cup_x z$. 
Let $\omega$ denote the homology class of $w$. 
We prove:

\f{\bf Claim 4.3.}    {\sl 
Consider $\omega\in$  $H_{2m+2}(V_2;\Bbb Z_2)$. 
The $\Bbb Z_2$-intersection number $\omega\cdot\omega$ in $V_2$ is one.}

{\bf Proof.} 
Let 
$W\x \{t\vert-1\leqq t\leqq 1\}$ 
be a tubular neighborhood of $W$ in $V_2$. 
The $(2n+1)$-cycle $x\x \{t=1\}$ in $S^{4m+3}_2$ 
bounds a $(2n+2)$-chain $\widetilde{z}$ in $S^{4m+3}_2$.
Let $w_t$ 
be  the cycle  

$(y\x \{t=1\})$ $\cup_{x\x \{t=1\}} \widetilde{z} $. 
\nl Let 
$\partial V_2\x \{s\vert 0\leqq s\leqq1\}$ 
be a collar neighborhood of $\partial V_2$ in $V_2$. 
Let $w_s$ 
be  the cycle 

$(y-[y\cap$ \{$\partial V_2\x$$ \{s\vert  0\leqq s\leqq1\}$]) 
 $\cup_{x\x \{s=1\}} (z\x\{s=1\}) $. 
\nl Then the followings hold.
(1) $w_s$ and  $w_t$ are in $\partial V_2$  
 and represent the same homology class $\omega$.
  (2)$w_s$ and  $w_t$  intersect transversely at odd points 
  because $\theta(x,x) $ is odd. 
  Therefore $\omega\cdot\omega$ is one.\qed

On the contrary to the above Claim 4.3, the following Claim 4.4 holds. 
This is absurdity. We prove Claim 4.4 and 
complete the proof of Proposition 4.1.

\f{\bf Claim 4.4.}    {\sl
The $\Bbb Z_2$-intersection number $\omega\cdot\omega$ is zero. }

{\bf Proof.} 
Since $V_2$ is a codimension one orientable submanifold in the parallelizable 
manifold $S^{4m+5}-$\{one point\} and $\partial V_2$ $\neq\phi$, 
$V_2$ is parallelizable. 
The fact that $\omega\cdot\omega$ is zero 
follows from the following elementary Lemma. 
This lemma is essentially same as  in P. 525 of [KM]. 
In fact, it is proved elementarily without using $Sq$-operators.

\f{\bf Lemma.}     {\sl (See e.g. P. 525 of [KM].)
Let $V$ be a compact parallelizable $2k$-manifold.
For an arbitrary $k$-homology class $\omega\in 
H_k(V ;\Bbb Z_2)$ 
the intersection number $\omega\cdot\omega$ is zero. }

\head 5.
A necessary condition for the realization of pair of $(4m+3)$-knots  
\endhead

In this section we prove: 
\proclaim{Proposition 5.1}
If a 4-tuple of ($4m+3, 4m+5$)-knots ($K_1, K_2,X_1,X_2$) 
 is realizable ($m\geqq0$), 
 then  $\sigma$($K_1$)=$\sigma$($K_2$).
\endproclaim

\f{\bf Proof of Proposition 5.1.}
Let $f:S^{4m+5}_1\coprod S^{4m+5}_2 \looparrowright S^{4m+7}$ 
be an immersion to realize ($K_1, K_2, X_1, X_2$). 
We abbreviate $X_i=f(S^{4m+5}_i)$ to $S^{4m+5}_i$.
There exist compact oriented ($4m+6$)-manifolds $V_1$ and $V_2$ such that 
$S^{4m+5}_i =\partial V_i \subset V_i \subset S^{4m+7}(i=1,2)$ and 
$V_1$ and $V_2$ intersect transversly.
Let $W$ be $V_1\cap V_2.$ Then 
$\partial W=\partial ( V_1\cap V_2)$
$=(\partial V_1\cap V_2)\cup (V_1\cap \partial V_2)$
$=(S^{4m+3}_1\cap V_2)\cup (V_1\cap S^{4m+3}_2).$
Let $F_1=S^{4m+3}_1\cap V_2$ and $F_2=V_1\cap S^{4m+3}_2$. 
Then $F_i$ in $S^{4m+5}_i$ is a Seifert hypersurface for $K_i$ ($i=1,2$).
Therefore 
$\sigma(K_1)-  \sigma(K_2)$
$=\sigma(K_1)+  \sigma(-K_2)$
$=\sigma(F_1)+  \sigma(-F_2)$
$=\sigma(\partial W)=0.$
The second equality holds by the definition of the signature. 
The third equality holds by Novikov aditivity. 
\qed

\head 5.A. The proof of Claim 2.2.2.
\endhead 

As we state under Theorem 2.2 in \S2, we prove Claim2.2.2 here. 
We complete the proof of Theorem 2.2.

Since $K_2$ and $K_3$ are cobordant,  
$\cases
\text{ 
 Arf$(K_2)$=Arf$(K_3)$ 
 }
 & \text{when $k$ is even}\\ 
 \text{ 
$\sigma (K_2)=\sigma (K_3)$ 
 }
 & \text{when $k$ is odd.} 
\endcases$
There exist 
$K'_1=K_1$, $K'_2$,..., $K'_q$, $K'_{q+1}=K_2$
and pass-move-charts $U_i$  of $K'_i$ $(i=1,...,q)$
and 
$K'_{i+1}$ is obtained from $K'_i$ 
by the high dimensional pass-move in $U_i$ $(i=1,...,q)$.
If the equality ($\dagger$) 
$\cases
\text{ 
 Arf$(K'_i)$=Arf$(K'_{i+1})$ 
 }
 & \text{when $k$ is even}\\ 
 \text{ 
$\sigma (K'_i)=\sigma (K'_{i+1})$ 
 }
 & \text{when $k$ is odd} 
\endcases$
holds for $i=1,...,q$, then the proof is completed.
By Lemma 3.4, the pair of ($2k+1$)-knots ( $K'_i$, $K'_{i+1}$ ) is realizable.
Therefore, by  Proposition 4.1 and 5.1, the equality ($\dagger$) holds.


\head 6.
A necessary and sufficient condition 
for the realization of 4-tuple of even dimensional knots  
\endhead

In this  section we prove Theorem 1.4.

We need the following Lemma 6.1.
\proclaim{Lemma 6.1}
Let $T$ be the trivial ($n+2$)-knot,  
$K'_2$ the trivial $n$-knot, and  $K_1$ a slice  $n$-knot.
($n\geqq1$). 
Then $(K_1, K_2', T, T)$ is realizable. 
\endproclaim

Before the proof of Lemma 6.1, we prove:

\f{\bf Claim.}{\sl 
If  Sublemma 3.5, Lemma 4.2 and 6.1 hold,  Theorem 1.4 holds. }

Proof.  
Let 
$K_2$ be a slice $n$-knot, 
$K'_1$ the trivial $n$-knot,
$X_i$ an arbitrary $(n+2)$-knot diffeomorphic to the standard $(n+2)$-sphere, 
and $T$ the trivial $(n+2)$-knot ($i=1,2$). 
By Lemma 6.1, $(K_1, K_2', T, T)$ 
and $(K_1', K_2, T, T)$ are realizable. 
By Lemma 4.2, 
$(K_1 \sharp K'_1,$ 
$ K'_2\sharp K_2,$ 
$ T   \sharp T,$ 
$ T   \sharp T$) 
=$(K_1, K_2, T, T)$ 
is realizable. 
By Sublemma 3.5, $(K_1, K_2, X_1, X_2 )$ is realizable.

We prove Lemma 6.1 to complete the proof of Theorem 1.4.

{\p Lemma 6.1.}
We define  $f:S^{n+2}_1\coprod S^{n+2}_2 \looparrowright S^{n+4}$ 
by using the $k$-twist spinning in \S 6 of [Z].
Prepare $D^{n,n-2}$ in  \S 6 of [Z],   
and  put $n$ there to be  $(n+3)$. 
As written there, 
regard 
$(S^{n+4}$, a $(n+2)$-knot) as 
($\partial D^{n+3,n+1}\times $$D^2$) $\cup$ 
($D^{n+3,n+1}\times$ $\partial D^2$).
Take $D^{n+3,n+1}$  as follows. 
Recall that (1)$D^{n+3,n+1}$ denote a set of  
$D^{n+3}$ and $D^{n+1}$ embedded in $D^{n+3}$,   
(2)
$D^{n+3}$ $\cap$ $D^{n+1}$ = $\partial D^{n+1}$    and 
$\partial D^{n+1}$ in $\partial D^{n+3}$ is the trivial $n$-knot. 
Regard $D^{n+3}$ as $D_s^{n+2}$ $\x[-1,1]$. 
Let  $D^{n+1}$ $\cap$ $\partial D^{n+3}$ $\subset$  ($D_s^{n+2}$ $\x\{-1\}$). 
Suppose that   
$D^{n+1}$ $\cap$ ($D_s^{n+2}$ $\x\{0\}$) 
in ($D_s^{n+2}$ $\x\{0\}$) defines $K_1$. 
Such  $D^{n+3,n+1}$ exists because $K_1$ is slice. 
Define $f\vert S^{n+2}_1$ 
so that $f(S^{n+2}_1)$ is the boundary of 
$D_s^{n+2}$ $\x[0,1]$ $\x$ $ \{\theta_0\}$, 
 where $\theta_0$ is a point in $\partial D^2$.
Define $f\vert S^{n+2}_2$ 
so that $f(S^{n+2}_2)$ coincides with what is made from 
$D^{n+1}$ by 1-twist-spinning.   
Then the following claim holds. we prove:

\f{\bf Claim.}
The immersion $f$ realizes the 4-tuple of ($n,n+2$)-knots $(K_1, K_2', T, T)$.

Proof. 
$f(S^{n+2}_1)$ is the boundary of the $(n+3)$-ball 
$D_s^{n+2}$ $\x[0,1]$ $\x$ $ \{\theta_0\}$. 
Therefore $f\vert S^{n+2}_1$ defines the trivial knot $T$. 
$f(S^{n+2}_1)$ is a 1-twist spun knot. 
By [Z] 1-twist spun knots are trivial. 
Therefore $f\vert S^{n+2}_2$ defines the trivial knot $T$. 
By the definition of the construction of $f$, 
the $n$-knot 
$f(S^{n+2}_1)$ $\cap$ $f(S^{n+2}_2)$ in $f(S^{n+2}_1)$ 
is 
$D^{n+1}$ $\cap$ ($D_s^{n+2}$ $\x\{0\}$) 
in ($D_s^{n+2}$. 
Therefore 
$f(S^{n+2}_1)$ $\cap$ $f(S^{n+2}_2)$ in $f(S^{n+2}_1)$ defines $K_1$. 
The $(n+1)$-disc 
$D^{n+1}$ $\cap$ ($D_s^{n+2}$ $\x[01]$) is called $D_l^{n+1}$. 
By the definition of the construction of $f$, 
$D_l^{n+1}$ is in $f(S^{n+2}_2)$. 
By the definition of the construction of $f$, 
the $n$-knot 
$f(S^{n+2}_1)$ $\cap$ $f(S^{n+2}_2)$ in $f(S^{n+2}_2)$ 
 is the boundary of $D_l^{n+1}$. 
 Therefore 
$f(S^{n+2}_1)$ $\cap$ $f(S^{n+2}_2)$ in $f(S^{n+2}_2)$ defines 
the trivial knot $K'_2$. 
Therfore $f$ is an immersion to realize  $(K_1, K_2', T, T)$. \qed


\np

\vskip3cm

Department of Mathematical Sciences, University of Tokyo

Komaba, Tokyo 153,   Japan
 
i33992\@ m-unix.cc.u-tokyo.ac.jp

\enddocument